\documentclass[12pt,a4paper]{article}
\usepackage{amsmath,mathrsfs,amsmath,amssymb,enumerate}
\usepackage{hyperref}
\usepackage{cleveref} 
\usepackage{dsfont}
\usepackage{comment}
\usepackage{indentfirst}
\usepackage{url}

\numberwithin{equation}{section}
\newtheorem{thm}{Theorem}[section]
\newtheorem{cor}[thm]{Corollary}
\newtheorem{lem}[thm]{Lemma}

\newtheorem{defn}[thm]{Definition}
\newtheorem{rem}[thm]{\bf Remark}

\mathchardef\lz="2D

\newcommand{\R}{\mathbb{R}^{N}}
\newcommand{\C}{C^{1}}
\newcommand{\CC}{C^{2}}

\newcommand{\HO}{{H}_0^{1}(\Omega)}

\newcommand{\LI}{{L}^{1}(\Omega)}
\newcommand{\LL}{L^{1}\lz solution}

\newcommand{\RR}{\mathbb{R}}

\allowdisplaybreaks

\crefname{equation}{Problem}{Problem}
\crefname{thm}{Theorem}{Theorem}
\crefname{lem}{Lemma}{Lemma}

\begin{document}

\baselineskip 22pt \setcounter{page}{1}
\title{\bf On Li-Lin's open problem\thanks{Supported by National Natural Science Foundation of China(No. 12371181).}
}
\author{{Zhi-Yun Tang, Xianhua Tang}\\
{\small \emph{School of Mathematics and Statistics, HNP-LAMA, Central South University,}}\\
{\small \emph{Changsha, Hunan 410083, People's Republic of China}}\\}
\date{}
\maketitle

{\bf Abstract}:
		In this paper, we give a first negative answer to a question
proposed by Li and Lin (Arch Ration Mech Anal 203(3): 943-968, 2012). Meanwhile we also give a second positive answer to the Li-Lin's open problem. The first positive answer was given by G. Cerami, X. Zhong and W. Zou (Calc. Var. Partial Differential Equations, 54(2): 1793-1829, 2015).



\section{Introduction}
Consider the Hardy-Sobolev's critical exponents problem with boundary singularities
		\begin{equation}
			\begin{cases}\label{eq1}
		-\Delta u=-\lambda |x|^{-s_{1}}|u|^{p-2}u+|x|^{-s_{2}}|u|^{q-2}u &\text { in } \Omega, \\
								u(x)=0 &\text { on } \partial \Omega,
			\end{cases}
		\end{equation}
where $\Omega \subset \R$, $N\geq3$, $\lambda\in \RR$, $0\leq s_1 < s_2 < 2$, $2^{*}(s)=\frac{2(N-s)}{N-2}$ for $0 \leq s \leq 2, 2< q\leq 2^{*}(s_2), 2< p\leq 2^{*}(s_1).$

In \cite{LL2012} Y. Li and C.-S. Lin asked an open question: Does \cref{eq1} have a positive solution when $\lambda>0$ and $p>q=2^{*}(s_2)$? In 2015, G. Cerami, X. Zhong and W. Zou provided some positive answers in \cite{CZZ2015}, which is the only answer so far. The main difficulty of the Li-Lin's open problem is the proof of the boundedness for the Palais-Smale sequence due to the nonlinearity terms do not satisfy the Ambrosetti-Rabinowitz type condition (see \cite{KT18EJDE} and the references here) except the one caused by the critical exponents. Motivated by the mentioned papers above, to the Li-Lin's open problem, we give a first negative answer by contradiction with the H\"{o}lder inequality, the Hardy inequality and the Young inequality, and a second positive answer by the method of sub-supersolutions. The main results are the following theorems.

\begin{thm}\label{thm1}
	Suppose that $\Omega \subset \R$ is a domain. Assume that $0\leq s_1 < s_2 < 2$, $p=2^{*}(s_1)$ and $q=2^{*}(s_2)$. Then there exists $\lambda_{1}>0$ such that \cref{eq1} has no nonzero solution for all $\lambda>\lambda_{1}$.

\end{thm}

\begin{rem}
For the Li-Lin's open problem, to our best knowledge, it seems to be the first work on the nonexistence of positive solution.
\end{rem}

\begin{rem}
The nonexistence of positive solution in Theorem \ref{thm1} does not need the smooth of $\partial \Omega$ and the curvature condition $H(0)<0$.
\end{rem}

\begin{thm}\label{thm2}
	Suppose that $\Omega \subset \R$ is a $\C$ bounded domain such that $0 \in \partial\Omega $. Assume that $\partial \Omega$ is $\CC$at $0$ and $H(0)<0$. Let $0\leq s_1 < s_2 < 2$, $q=2^{*}(s_2)$, $2^{*}(s_1)-\frac{N(s_2-s_1)}{(N-2)(N+1-s_2)}< p\leq 2^{*}(s_1)$ and $2^{*}(s_1)-\frac{s_2-s_1}{N-2}\leq p$. Let $\lambda_{*}=\sup \{\lambda \in \RR\ |$\ \cref{eq1}\ has a positive solution\}. Then $\lambda_{*}>0$ and \cref{eq1} has at least a positive solution for all $\lambda \in (0,\lambda_{*})$.
\end{thm}

\begin{cor}\label{cor}
	Suppose that $\Omega \subset \R$ is a $\C$ bounded domain such that $0 \in \partial\Omega $. Assume that $\partial \Omega$ is $\CC$at $0$ and $H(0)<0$. Let $0\leq s_1 < s_2 < 2$, $q=2^{*}(s_2)$ and $p= 2^{*}(s_1)$. Then there exists $\lambda_{*}>0$ such that \cref{eq1} has no positive solution for all $\lambda>\lambda_{*}$ and \cref{eq1} has at least a positive solution for all $\lambda \in (0,\lambda_{*})$.
\end{cor}

\begin{rem}
Theorem \ref{thm2} is the first result on the existence of positive solution by the method of sub-supersolutions, which had been obtained by G. Cerami, X. Zhong and W. Zou in \cite{CZZ2015} for $\lambda>0$ sufficient small by the variational method.
\end{rem}

\begin{rem}
We guess that the $\lambda_{*}$ in Theorem \ref{thm2} is infinite when $p<2^{*}(s_1)$.
\end{rem}

\section{Preliminaries}
First, we recall positive answers for Li-Lin's open problem in \cite{CZZ2015}, which will be used in the late.
\begin{lem}\label[lem]{lem1}
	\emph{(see \cite{CZZ2015}, Theorem 1.5)}
	Suppose that $\Omega \subset \R$ is a $\C$ bounded domain such that $0 \in \partial\Omega$. Assume that $\partial \Omega$ is $\CC$ at $0$ and $H(0)<0$. Let $0\leq s_1 < s_2 < 2$ and $2^{*}(s_2)<p\leq 2^{*}(s_1)$. Then there exists $\lambda_{0}>0$ such that
	\cref{eq1} has at least a positive solution for all $\lambda \in (0,\lambda_{0})$.
\end{lem}

Then, we recall the sub-supersolution method in \cite{MP2008}.

\begin{defn}\label{def1}
	\emph{(see \cite{MP2008}, P2430, Definition 1.1)}
	A function u is an $\LL$ of
	\begin{equation}\label{eq2}
		\begin{cases}
			-\Delta u=f(x,u) &\text { in } \Omega, \\
			u(x)=0 &\text { on } \partial \Omega,
		\end{cases}
	\end{equation}
	where $\Omega \subset \R$ is a smooth bounded domain and $f:\Omega \times \mathbb{R} \rightarrow \mathbb{R}$ is a Carath$\acute{e}$odory
	function, if
\begin{flalign*}
		(i)\ &u \in \LI;\ &\\
		(ii)\ &f(\cdotp,u)\rho_0 \in\LI;\\
		(iii)\ &
\end{flalign*}
\begin{equation}\label{eq22}
	-\int_{\Omega}u\Delta\varphi dx
	=\int_{\Omega}f(x,u)\varphi dx \ \ \forall \varphi \in C_0^2(\overline{\Omega}).
\end{equation}
Here, $\rho_0(x)=d(x,\partial\Omega),\ \forall x \in \Omega,\  and \  C_0^2(\overline{\Omega})=\{ \varphi \in C^2(\overline{\Omega});\ \varphi=0\ on\ \partial \Omega \}$.
\end{defn}

We also consider $L^1 \lz sub$ and $L^1\lz super solutions$ in analogy with this definition. For instance, $u$ is an $L^1 \lz sub solutions$ of \cref{eq2} if $u$ satisfies (i)-(iii) with ``$\leq$'' instead of ``$=$'' in (\ref{eq22}). A similar convention for sub and supersolutions.

\begin{lem}\label[lem]{sxj1.1}
	\emph{(see \cite{MP2008}, P2430, Theorem 1.1)}
	Let $v_1$ and $v_2$ be a sub and a supersolution of \cref{eq2}, respectively. Assume that $v_1 \leq v_2$ a.e. and
	\begin{equation}\label{eq23}
		f(\cdotp,v)\rho_0 \in\LI\ \ \ for\ every\ v\in \LI\ such\ that\ v_1 \leq v \leq v_2\ a.e.
	\end{equation}
	Then, there exists a solution $u$ of \cref{eq2} in $[v_1,v_2]$.
\end{lem}

Now we give some lemmas which will be needed in the proofs of our main results.

\begin{lem}\label[lem]{alpha}
Suppose that $\Omega \subset \R$ is a bounded domain such that $0 \in \overline{\Omega} $. Let $v(x)=|x|^{-\alpha}$, $x \in \overline{\Omega}$. Then $v \in\LI $ for $\alpha <N$ and $v \in L^1(\partial \Omega) $ for $\alpha < N-1$.
\end{lem}

{\bf Proof}\ \ \ \
First, we prove that $v \in\LI $. In fact, due to $N-\alpha> 0$, one has
\begin{align*}
	\int_{\Omega}|x|^{-\alpha}\ dx&\leq \int_{B_{1}(0)}|x|^{-\alpha} dx+\int_{\Omega\setminus B_{1}(0)}|x|^{-\alpha} dx\\
	&\leq \omega_N \int_0^{1}r^{-\alpha }r^{N-1}dr+\int_{\Omega\setminus B_{1}(0)}|x|^{-\alpha} dx\\
	&= \frac{\omega_N}{N-\alpha}+\int_{\Omega\setminus B_{1}(0)}|x|^{-\alpha} dx\\
	&< +\infty,
\end{align*}
where $\omega_N$ is the volume of unit ball in $\R$, hence $v \in\LI $.

Next we prove that $v \in L^1(\partial \Omega) $.  Because that $|\nabla v(x)|=|\alpha||x|^{-\alpha-1}$, we have $v \in W^{1,1}(\Omega)$ for $\alpha < N-1$. It follows from Sobolev trace theorem that $W^{1,1}(\Omega) \hookrightarrow L^{1}(\partial\Omega)$, which completes our proving.  $\hfill\Box$

\begin{lem}\label[lem]{super}
Suppose that $\Omega \subset \R$ is a bounded domain such that $0 \in \partial \Omega $. Let $\alpha\in [0, N-2]$ and $v(x)=|x|^{-\alpha}$, $x \in \overline{\Omega}$. Then $-\int_{\Omega}v\Delta\varphi dx
\geq 0$ for any non-negative function $\varphi \in C_0^2(\overline{\Omega})$.
\end{lem}

{\bf Proof}\ \ \ \ Choose $ \psi \in C^2_0(\mathbb{R})$ such that $\psi (s)=1$ for $|s|\leq 1$ and $\psi(s)=0$ for $|s|\geq 2$. On one hand, we have
\begin{align}\label{2.5}
	\lim_{\varepsilon \to 0}\left(-\int_{\Omega}\psi \left( \varepsilon^{-1}|x| \right) |x|^{-\alpha}\Delta \varphi dx\right)=0
\end{align}
by Lebesgue Dominated Convergence Theorem and \cref{alpha}. On the other hand, we obtain
\begin{align*}
	&-\int_{\Omega}\left(1-\psi\left( \varepsilon^{-1}|x| \right) \right) |x|^{-\alpha}\Delta \varphi dx\\
	=&-\int_{\Omega}\varphi\Delta \left[\left(1-\psi\left( \varepsilon^{-1}|x| \right) \right) |x|^{-\alpha} \right]dx\\
	&+\int_{\partial\Omega} \left\{\varphi\frac{\partial}{\partial \nu} \left[\left(1-\psi\left( \varepsilon^{-1}|x| \right) \right) |x|^{-\alpha} \right]-\left(1-\psi\left( \varepsilon^{-1}|x| \right)\right)|x|^{-\alpha}\frac{\partial \varphi }{\partial \nu }\right\}ds\\
	=& -\int_{\Omega}\varphi\Delta\left[\left(1-\psi\left( \varepsilon^{-1}|x| \right) \right) |x|^{-\alpha} \right]dx-\int_{\partial\Omega} \left(1-\psi\left( \varepsilon^{-1}|x| \right)\right)|x|^{-\alpha}\frac{\partial \varphi }{\partial \nu }ds
\end{align*}
by Green's second identity. Define $$f_\varepsilon (x)=-\Delta\left(\left(1-\psi\left( \varepsilon^{-1}|x| \right) \right) |x|^{-\alpha}\right)$$ and $f(x)=\alpha(N-2-\alpha)|x|^{-\alpha-2}$. Then one has
\begin{align*}
f_\varepsilon (x)=&-div\left(\nabla\left(\left(1-\psi\left( \varepsilon^{-1}|x| \right) \right) |x|^{-\alpha}\right)\right) \\
	=&\varepsilon^{-1}div\left(|x|^{-\alpha-1}\psi'\left( \varepsilon^{-1}|x| \right)x\right)+\alpha  div\left(\left(1-\psi\left( \varepsilon^{-1}|x| \right) \right)|x|^{-\alpha-2} x\right)\\
	=&\varepsilon^{-2}|x|^{-\alpha}\psi''\left( \varepsilon^{-1}|x| \right)+(N-1-2\alpha)\varepsilon^{-1}|x|^{-\alpha-1}\psi'\left( \varepsilon^{-1}|x| \right)\\ &+\alpha(N-2-\alpha)\left(1-\psi\left( \varepsilon^{-1}|x|\right) \right) |x|^{-\alpha-2}
\end{align*}
Let $$h(x)=\left(4\|\psi''\|_\infty+2(N-1+2\alpha)\|\psi'\|_\infty+\alpha(N-2-\alpha)\left(1+\|\psi\|_\infty \right)\right) |x|^{-\alpha-2}.$$ Then one has $h\varphi\in \LI$, $|f_\varepsilon(x)|\varphi\leq h(x)\varphi$ for $\varepsilon \in (0,1]$ and $x \in \Omega$, and $f_\varepsilon(x)\varphi\to f(x)\varphi$ as $\varepsilon \to 0$ for $x\in \Omega$, which implies that
\begin{align*}
	&-\lim_{\varepsilon \to 0}\int_{\Omega}\left(1-\psi\left( \varepsilon^{-1}|x| \right) \right) |x|^{-\alpha}\Delta \varphi dx\\
&=\lim_{\varepsilon \to 0}\int_{\Omega}f_{\varepsilon}(x)\varphi dx-\lim_{\varepsilon \to 0}\int_{\partial\Omega}\left(1-\psi\left( \varepsilon^{-1}|x| \right) \right)|x|^{-\alpha}\frac{\partial \varphi}{\partial\nu}ds\\
	&=\int_{\Omega}f(x)\varphi dx-\int_{\partial\Omega}|x|^{-\alpha}\frac{\partial \varphi}{\partial\nu}ds
\end{align*}
by Lebesgue Dominated Convergence Theorem, where we use \cref{alpha} due to $\alpha<N-1$. Hence, one obtains
\begin{align*}
	-\int |x|^{-\alpha}\Delta \varphi dx=&-\lim_{\varepsilon \to 0}\int_{\Omega}\psi_{\varepsilon}(x) |x|^{-\alpha}\Delta \varphi dx-\lim_{\varepsilon \to 0}\int_{\Omega}\left(1-\psi_{\varepsilon}(x)\right) |x|^{-\alpha}\Delta \varphi dx\\
	=&\alpha(N-2-\alpha)\int_{\Omega}|x|^{-\alpha-2}\varphi dx-\int_{\partial\Omega}|x|^{-\alpha}\frac{\partial \varphi}{\partial\nu}ds\\
	\geq&0
\end{align*}
by \eqref{2.5} and the fact that $\frac{\partial \varphi}{\partial\nu}|_{\partial\Omega}\leq 0$ for $\varphi(x)\geq 0$ on $\overline{\Omega}$, which complete the proof.
$\hfill\Box$
\\

The next regularity result is inspired by \cite{S1964,CC1993}.
\begin{lem}\label[lem]{regularity}
	Suppose that $\Omega \subset \R$ is a bounded domain such that $0 \in \partial\Omega $, $\mu>0$, $0\leq s_1 < s_2 <2$, $q<p\leq 2^{*}(s_1)$ and $2<q\leq 2^{*}(s_2)$.  Assume that $w\in \HO$ satisfies
	\begin{align}\label{ce1}
	-\Delta w=-\mu |x|^{s_{1}}|w|^{p-2}w+|x|^{-s_{2}}|w|^{q-2}w &\text { in } \Omega.
\end{align}
	Then $w\in C(\overline{\Omega}) \cap C^{2}(\Omega)$.
\end{lem}

{\bf Proof}\ \ \ \ It follows from the Moser iteration technique and the regularity theory of elliptic equations that $w \in C^{2}(\Omega)$.

Choose
$\eta \in C_{0}^{\infty}\left(B_{\rho}(0),[0,1]\right)$
with $\eta=1$ in $B_{\rho_0}(0)$ and $|\nabla \eta| \leq \frac{4}{\rho-\rho_0}$ ($ 0<\rho_0<\rho<1$). For each $k > 1$, we define
\begin{align*}
	w_k(x)=
	\begin{cases}
		w(x)\ \ \ \ \  &\mathrm{if}~~|w(x)|\leq k,\\
		\pm k   &\mathrm{if}~~\pm w(x)>k.
	\end{cases}
\end{align*}
For $\beta > 1$, we use $\varphi_k = \eta^2 |w_k|^{2(\beta-1)}w$ as a test function in \eqref{ce1} to obtain
\begin{align}\label{4.2}
	&\int_{\Omega}
	\left(2\eta |w_k|^{2(\beta-1)}w \nabla w\cdot \nabla\eta +\eta^2 |w_k|^{2(\beta-1)} |\nabla w|^2+2(\beta-1)\eta^2 |w_k|^{2(\beta-1)}|\nabla  w_k|^2\right)dx \nonumber\\
	=&\int_{\Omega} \nabla w \cdot\nabla \varphi_kdx \nonumber\\
	=&\int_{\Omega}\eta^2|w_k|^{2(\beta-1)} \left(-\mu |x|^{-s_{1}}|w|^{p}+|x|^{-s_{2}}|w|^q \right)dx.
\end{align}
By Young's inequality we have
\begin{align}\label{3.16}
	\left|2 \int_{\Omega}  \eta |w_k|^{2(\beta-1)}w \nabla w\cdot \nabla\eta dx \right|
	\leq & \frac{1}{2} \int_{\Omega} \eta^{2} |w_{k}|^{2(\beta-1)}|\nabla w|^{2} d x+2\int_{\Omega}w^{2} |w_{k}|^{2(\beta-1)}|\nabla \eta|^{2} d x.
\end{align}
Then \eqref{3.16} is substituted in \eqref{4.2},
\begin{align}\label{2.7}
	&\frac{1}{2}\int_{\Omega}
	\left(\eta^2 |w_k|^{2(\beta-1)} |\nabla w|^2+2(\beta-1)\eta^2 |w_k|^{2(\beta-1)}|\nabla  w_k|^2\right)dx\nonumber\\
	\leq& 2\int_{\Omega}w^{2} |w_{k}|^{2(\beta-1)}|\nabla \eta|^{2} d x+\int_{\Omega}  \eta^2|w_k|^{2(\beta-1)}|x|^{-s_{2}}|w|^q dx.
\end{align}
Now we recall the Hardy-Sobolev's inequality,
\begin{align}\label{ckn}
	\left(\int_{\Omega}|x|^{-s}|v|^q d x\right)^{\frac{2}{q}} \leq C_{s,q}
	\int_{\Omega}|\nabla v|^{2} d x \ \ \ \ \mathrm{ for~all}~ v \in H_0^{1}\left(\Omega\right),
\end{align}
where $s\in [0, 2], 2\leq q\leq 2^*(s)=\frac{2 (N-s)}{N-2}$ and $C_{s,q}$ is a positive
constant depending on $s, q.$

Choosing $v=\eta |w_{k}|^{\beta-1}w $ in \eqref{ckn}, together with \eqref{2.7}, we derive
\begin{align}\label{2.9}
	&\int_{\Omega} \left|\nabla v\right|^{2}dx\nonumber\\&=\int_{\Omega} \left|\nabla (\eta |w_{k}|^{\beta-1}w)\right|^{2}dx\nonumber\\
	&\leq  \int_{\Omega}\left| |w_{k}|^{\beta-1}w\nabla \eta+
	\eta |w_{k}|^{\beta-1}\nabla w+ (\beta-1)\eta |w_{k}|^{\beta-1}\nabla w_k\right|^{2}dx\nonumber\\
	&\leq  3\int_{\Omega} \left( |w_{k}|^{2(\beta-1)}w^2|\nabla \eta|^2+
	\eta^2 |w_{k}|^{2(\beta-1)}|\nabla w|^2+ (\beta-1)^2\eta^2|w_{k}|^{2(\beta-1)}|\nabla w_k|^2\right)dx\nonumber\\
	&\leq  24\beta\int_{\Omega}\left(w^{2} |w_{k}|^{2(\beta-1)}|\nabla \eta|^{2}+  \eta^2|w_k|^{2(\beta-1)}|x|^{-s_{2}}|w|^q \right)dx.
\end{align}
Therefore, by the absolute continuity of the integral, there exists a small $\rho>0$ such that
\begin{align}\label{2.10}
	&24\beta\int_{\Omega} \eta^2 |w_k|^{2(\beta-1)}|x|^{-s_{2}}|w|^q dx\nonumber\\=&
	24\beta\int_{B_{\rho}(0)}  |x|^{-s_{2}}|w|^{q-2}  v^2dx\nonumber\\
	\leq& 24\beta\left(\int_{B_{\rho}(0)}|x|^{-s_{2}} |w|^q dx\right)^{\frac{q-2}{q}}\left(\int_{\Omega} |x|^{-s_{2}} |v|^{q} dx\right)^{\frac{2}{q}}\nonumber\\
	\leq& \frac{1}{2}  \int_{\Omega} \left|\nabla v\right|^{2}dx.
\end{align}
It follows from \eqref{2.9} and \eqref{2.10} that
\begin{align*}
	\left(\int_{\Omega}
	\left|\eta |w_{k}|^{\beta-1}w\right|^{2^*}dx\right)^{\frac{2}{2^*}}\leq  C_{0}\int_{\Omega} \left|\nabla v\right|^{2}dx \leq 48 \beta C_{0}\int_{\Omega} w^{2} |w_{k}|^{2(\beta-1)}|\nabla \eta|^{2}dx.
\end{align*}
Using the choice of the cut-off function $\eta$, we deduce
\begin{align*}
	\left(\int_{B_{\rho_0}(0)\cap \Omega}
	\left| |w_{k}|^{\beta-1}w\right|^{2^*}dx\right)^{\frac{2}{2^*}}
	\leq \frac{768 \beta C_{0}}{(\rho-\rho_0)^2} \int_{B_{\rho}(0)\cap \Omega}w^{2} |w_{k}|^{2(\beta-1)} dx.
\end{align*}
Letting $k\to\infty$, we
have
\begin{align}\label{3.4}
	\left(\int_{B_{\rho_0}(0)\cap \Omega}
	|w|^{2^* \beta} dx\right)^{\frac{1}{2^*}}\leq \frac{16\sqrt{3 \beta C_{0}}}{\rho-\rho_0} \left(\int_{B_{\rho}(0)\cap \Omega}|w|^{2\beta} d x\right)^{\frac{1}{2}}.
\end{align}
By the Moser iteration technique and the regularity theory of elliptic equations, we have $w \in C(\overline{\Omega}\setminus B_{\rho_0}(0)) $, which implies that $w\in L^{\frac{(2^*)^2}{2}}(\Omega)$ by (\ref{3.4}) with $\beta=\frac{2^*}{2}$.

Now we prove that $w\in L^r(\Omega)$ for all $r>1$. Let $a(x)=-\mu |x|^{-s_1}|w|^{p-2}+|x|^{-s_2} |w|^{q-2}$. Then $a\in L^{\frac{N}{2}}(\Omega)$. In fact,
let $q_1=\frac{(2^*)^2}{N(q-2)}$ and $q_2=\frac{(2^*)^2}{(2^*)^2-N(q-2)}$, then $-\frac{N}{2}s_2q_2>-N$, which implies that
\begin{align*}
	\int_{\Omega}{|x|^{-\frac{N}{2}s_2}}{|w|^{\frac{N}{2}(q-2)}} dx &\leq \left(\int_{\Omega}|x|^{-\frac{N}{2}s_2q_2}dx\right)^{\frac{1}{q_2}}\left(\int_{\Omega}|w|^{\frac{N}{2}(q-2){q_1}}dx\right)^{\frac{1}{q_1}}\\
	&=\left(\int_{\Omega}|x|^{-\frac{N}{2}s_2q_2}dx\right)^{\frac{1}{q_2}}\left(\int_{\Omega}|w|^\frac{(2^*)^2}{2}dx\right)^{\frac{1}{q_1}}\\
	&< +\infty,
\end{align*}
that is, ${|x|^{-s_2}}{|w|^{q-2}} \in L^{\frac{N}{2}}(\Omega)$. In a similar way, one has ${|x|^{-s_1}}{|w|^{p-2}} \in L^{\frac{N}{2}}(\Omega)$. Hence $a\in L^{\frac{N}{2}}(\Omega)$. It follows from Brezis-Kato Theorem (see \cite{BK1979}), $w\in L^r(\Omega)$ for all $r>1$.

For $r\in (\frac{N}{2}, \frac{N}{s_2})$, choose $q_3>1$ such that $-rs_2q_3>-N$. Let $q_4=\frac{q_3}{q_3-1}$. Then one has
\begin{align*}
	\int_{\Omega}\left(|x|^{-s_2}{|w|^{q-1}}\right)^r dx &\leq \left(\int_{\Omega}|x|^{-rs_2q_3}dx\right)^{\frac{1}{q_3}}\left(\int_{\Omega}|w|^{r(q-1){q_4}}dx\right)^{\frac{1}{q_4}} < +\infty
\end{align*}
by \cref{alpha}. Let $b(x)=-\mu {|x|^{-s_1}}{|w|^{p-2}w} + {|x|^{-s_2} }{|w|^{q-2}w}$. Then $b\in L^r(\Omega)$ for some $r>\frac{N}{2}$.
It follows from the regularity theory of elliptic equations that $w \in W^{2,r}(\Omega)$ for some $r>\frac{N}{2}$, which implies that $w\in C(\overline{\Omega})$ by Sobolev's imbedding theorem. This completes the proof.
$\hfill\Box$

\section{Proofs of Theorems \ref{thm1} and \ref{thm2}}

\subsection{Proofs of Theorem \ref{thm1}}

Suppose that $u$ is a nonzero solution of \cref{eq1}. Then $I'(u)=0$, where
\begin{equation*}
	I(u)=\dfrac{1}{2}\int_{\Omega}|\nabla u|^{2}dx
	+\frac{\lambda}{2^{*}(s_{1})} \int_{\Omega}|x|^{-s_{1}}|u|^{2^{*}(s_{1})}dx
	-\frac{1}{2^{*}(s_{2})}\int_{\Omega}|x|^{-s_{2}}|u|^{2^{*}(s_{2})} dx
\end{equation*}
on $H^1_0(\Omega)$. Hence $\langle I'(u),u\rangle=0$, that is,
\begin{equation}\label{2.1}
	\int_{\Omega}|\nabla u|^{2}dx
	+\lambda\int_{\Omega}|x|^{-s_{1}}|u|^{2^{*}(s_{1})}dx
	=\int_{\Omega}|x|^{-s_{2}}|u|^{2^{*}(s_{2})} dx.
\end{equation}
By the H\"{o}lder inequality and the Hardy inequality, we have
\begin{align*}
\int_{\Omega}|x|^{-s_{2}}|u|^{2^{*}(s_{2})} dx&=	\int_{\Omega} |x|^{{-\tfrac{s_{1}(2-s_2)}{2-s_1}}}|u|^{\tfrac{2^{*}(s_{1})(2-s_2)}{2-s_1}}   \cdot |x|^{-{\tfrac{2(s_2-s_1)}{2-s_1}}}|u|^{{\tfrac{2(s_2-s_1)}{2-s_1}}}dx\\
&\leq\left(\int_{\Omega}|x|^{-s_{1}}|u|^{2^{*}(s_{1})} dx\right)^{\tfrac{2-s_2}{2-s_1}}
\left(\int_{\Omega}|x|^{-2}|u|^{2}dx \right)^{\tfrac{s_2-s_1}{2-s_1}}\\
&\leq	\left(\int_{\Omega}|x|^{-s_{1}}|u|^{2^{*}(s_{1})} dx\right)^{\tfrac{2-s_2}{2-s_1}}
	\left( C_N \int_{\Omega}|\nabla u|^2dx \right)^{\tfrac{s_2-s_1}{2-s_1}},
\end{align*}
where $C_N=\left(\dfrac{2}{N-2}\right)^2$. It follows from  the Young inequality that
\begin{align*}
	&\int_{\Omega}|x|^{-s_{2}}|u|^{2^{*}(s_{2})} dx\\&\leq	 \left(\left(C_N^{-1}\frac{2-s_1}{s_2-s_1}\right)^{\tfrac{s_1-s_2}{2-s_1}}\left(\int_{\Omega}|x|^{-s_{1}}|u|^{2^{*}(s_{1})} dx\right)^{\tfrac{2-s_2}{2-s_1}}\right)
	\left(\left(\dfrac{2-s_1}{s_2-s_1}\right)\int_{\Omega}|\nabla u|^2dx \right)^{\tfrac{s_2-s_1}{2-s_1}}\\
	&\leq \left({\dfrac{2-s_2}{2-s_1}}\right)\left(C_N^{-1}\frac{2-s_1}{s_2-s_1}\right)^{\tfrac{s_1-s_2}{2-s_2}}	\int_{\Omega}|x|^{-s_{1}}|u|^{2^{*}(s_{1})} dx+\int_{\Omega}|\nabla u|^2dx,
\end{align*}
which implies that
\begin{align*}
	\lambda&\leq \lambda_1\stackrel{\triangle}{=}\left({\dfrac{2-s_2}{2-s_1}}\right)\left(C_N^{-1}\tfrac{2-s_1}{s_2-s_1}\right)^{\tfrac{s_1-s_2}{2-s_2}}
\end{align*}
by (\ref{2.1}). Hence the theorem holds.
$\hfill\Box$

\subsection{Proof of Theorem \ref{thm2}}

In this subsection, we will give the proof of \cref{thm2} by the method of sub-supersolutions. In fact, Theorem \ref{thm2} is a corollary of the following theorem with $q=2^{*}(s_2)$.

\begin{thm}\label{thm3}
	Suppose that $\Omega \subset \R$ is a bounded domain such that $0 \in \partial\Omega $. Let $0\leq s_1 < s_2 < 2$, $p\leq 2^{*}(s_1)$, $q\leq 2^{*}(s_2)$, $2<q<\frac{N+1-s_2}{N+1-s_1}p+\frac{s_2-s_1}{N+1-s_1}$, $q\leq p-\frac{s_2-s_1}{N-2}$ and $\lambda_{*}=\sup \{\lambda \in \RR\ |$\ \cref{eq1}\ has a positive solution\}.  Assume that $\lambda_{*}>0$. Then \cref{eq1} has at least a positive solution for all $\lambda \in (0,\lambda_{*})$.
\end{thm}

{\bf Proof}\ \ \ \ By the definition of $\lambda_{*}$, for every $\lambda \in (0,\lambda_{*})$, there exists $\mu \in (\lambda,\lambda_{*})$ such that \cref{eq1} with $\lambda=\mu$ has a positive solution $u_{\mu}$.

Let $v_1=u_{\mu}$. Then we prove that $v_1$ is a subsolution of \cref{eq1}. In fact, it follows from \cref{regularity} that $u_\mu \in C(\overline{\Omega})\subset\LI$. It is obvious that $g_{\lambda}(x,u_{\mu})\rho_0\in \LI$, where
$$g_\lambda(x,t)=-\lambda{|x|^{-s_{1}}}{|t|^{p-2}t} +{|x|^{-s_{2}}}{|t|^{q-2}t} $$
for $x\in \Omega$ and $t\in \RR$. Moreover, one obtains $u_\mu \in C_0(\overline{\Omega}) \cap C^2(\Omega)$ from \cref{regularity}. Then by Green's second identity, we have
\begin{eqnarray*}-\int_{\Omega}u_{\mu}\Delta\varphi dx&=&-\int_{\Omega}\varphi\Delta u_{\mu} dx+\int_{\partial\Omega}\left(\varphi\frac{\partial u_\mu}{\partial \nu}-u_\mu\frac{\partial \varphi}{\partial \nu}\right)ds\\
&=&\int_{\Omega}\varphi g_{\mu}(x,u_{\mu})dx\\
&\leq& \int_{\Omega}\varphi g_{\lambda}(x,u_{\mu})dx
\end{eqnarray*}
for all $\varphi \in C_0^2(\overline{\Omega})$ with $\varphi(x)\geq 0$ on $\overline{\Omega}$. Hence $v_1$ is a subsolution of \cref{eq1}.

Next we prove that $v_2(x)\stackrel{\vartriangle }{=}M|x|^{-\alpha}$ is a supersolution of \cref{eq1}, where $\alpha\stackrel{\vartriangle }{=}\frac{s_2-s_1}{p-q}\leq N-2$, which follows from $q\leq p-\frac{s_2-s_1}{N-2}$, and $$M=\max\left\{ \lambda^{-\frac{1}{p-q}},||u_\mu||_{\infty}\cdot\sup \left\{ |x|^{\alpha} \big| x \in \overline{\Omega} \right\}\right\}.$$
In fact,
we have $v_2\in \LI$ and
\begin{align*}
	g_\lambda (x,v_2(x))=&-\lambda M^{p}|x|^{-\alpha(p-1)}|x|^{-s_1}+M^{q}|x|^{-\alpha(q-1)}|x|^{-s_2}\\
	=&-M^{q}|x|^{-\alpha(p-1)}|x|^{-s_1}(\lambda M^{p-q}-|x|^{\alpha(p-q)-s_2+s_1})\\
	=&-M^{q}|x|^{-\alpha(p-1)}|x|^{-s_1}(\lambda M^{p-q}-1)\\ \leq& 0,
\end{align*}
by \cref{super}, one has
$$-\int_{\Omega}v_2\Delta\varphi dx \geq 0 \geq \int_{\Omega}\varphi g_{\lambda}(x,v_2)dx$$
for all $\varphi \in C_0^2(\overline{\Omega})$ with $\varphi(x)\geq 0$ on $\overline{\Omega}$. Note $g_\lambda (x,v_2)\rho_0 \in \LI$ by $\alpha(p-1)+s_1-1<N$ which follows from $q<\frac{N+1-s_2}{N+1-s_1}p+\frac{s_2-s_1}{N+1-s_1}$. Hence $v_2$ is a supersolution of \cref{eq1}.

Moreover one has
$$
v_2(x)=M|x|^{-\alpha}\geq M\left\{\sup \left\{ |x|^{\alpha} \big| x \in \overline{\Omega} \right\}\right\}^{-1}\geq ||u_{\mu}||_\infty \geq u_\mu(x) =v_1(x)
$$
for all $x\in \overline{\Omega}$. It is obvious that (\ref{eq23}) holds. Hence the proof is complete by \cref{sxj1.1}.
$\hfill\Box$

{\bf Proof of Theorem \ref{thm2}}\ \ \ \ Note that $q=2^{*}(s_2)$. By Lemma \ref{lem1}, we have $\lambda_{*}>0$. Now Theorem \ref{thm2} follows from Theorem \ref{thm3}.
$\hfill\Box$

\bibliographystyle{unsrt}
\bibliography{ref}

\end{document}